\documentclass{IEEEtran4PSCC}

\ifCLASSINFOpdf
   \usepackage[pdftex]{graphicx}
\else
   \usepackage[dvips]{graphicx}
\fi

\usepackage{array}
\usepackage{booktabs}
\usepackage{cite}

\usepackage[colorlinks=true, linkcolor=blue, urlcolor=blue, citecolor=blue]{hyperref}

\usepackage{url}
\usepackage[cmex10]{amsmath}
\usepackage{amsfonts,amsmath,amssymb,amsthm}
\usepackage{xcolor}

\newcommand{\DS}{\displaystyle}
\DeclareMathOperator*{\st}{s.t.}

\makeatletter
\let\old@ps@headings\ps@headings
\let\old@ps@IEEEtitlepagestyle\ps@IEEEtitlepagestyle
\def\psccfooter#1{%
    \def\ps@headings{%
        \old@ps@headings%
        \def\@oddfoot{\strut\hfill#1\hfill\strut}%
        \def\@evenfoot{\strut\hfill#1\hfill\strut}%
    }%
    \def\ps@IEEEtitlepagestyle{%
        \old@ps@IEEEtitlepagestyle%
        \def\@oddfoot{\strut\hfill#1\hfill\strut}%
        \def\@evenfoot{\strut\hfill#1\hfill\strut}%
    }%
    \ps@headings%
}
\makeatother

\begin{document}

\title{Hedging against Black Swans in Day-Ahead\\Energy Markets}

\author{
\IEEEauthorblockN{Liviu Aolaritei$^{*, \dagger,1}$, Boubacar Bangoura$^{*,2}$, Saverio Bolognani$^2$, Nicolas Lanzetti$^2$, and Florian D\"orfler$^2$}

\IEEEauthorblockA{$^1$Department of EECS, University of California, Berkeley, USA}
\IEEEauthorblockA{$^2$Automatic Control Laboratory, ETH Zurich, Switzerland}
}

\maketitle

\begin{abstract}
Renewable generators must commit to day-ahead market bids despite uncertainty in both production and real-time prices. While forecasts provide valuable guidance, rare and unpredictable extreme events\,---\,so-called black swans\,---\,can cause substantial financial losses. This paper models the nomination problem as an instance of optimal transport-based distributionally robust optimization (OT-DRO), a principled framework that balances risk and performance by accounting not only for the severity of deviations but also for their likelihood. The resulting formulation yields a tractable, data-driven strategy that remains competitive under normal conditions while providing effective protection against extreme price spikes. Using four years of Finnish wind farm and market data, we demonstrate that OT-DRO consistently outperforms forecast-based nominations and significantly mitigates losses during black swan events.
\end{abstract}


\thanksto{\noindent $^*:$ Equal contribution. \\
$^\dagger:$ Corresponding author ({\tt liviu.aolaritei@berkeley.edu}).
\\\noindent Liviu Aolaritei acknowledges support from the Swiss National Science Foundation through the Postdoc.Mobility Fellowship (grant agreement P500PT\_222215). This work was supported as a part of NCCR Automation, a National Centre of Competence in Research, funded by the Swiss National Science Foundation (grant number 51NF40\_225155).}

\section{Introduction}
\label{sec:intro}

Electricity markets are undergoing a profound transformation as renewable generation drives the global energy transition. Among renewable sources, wind power has emerged as a cornerstone of decarbonization, offering large-scale, low-cost, and carbon-free energy. Policy targets and evolving market designs are accelerating its integration into power systems across many countries. However, the very characteristics that make wind appealing also make it challenging to trade: its production depends entirely on the weather, which is inherently uncertain and exhibits chaotic variability even over short horizons.

In day-ahead electricity markets, wind producers must submit their bids by noon on day $D$, committing to a generation schedule for all 24 hours of day $D+1$. This nomination must be made without knowing either future production or market prices. Revenues are settled after delivery: nominated energy is remunerated at the day-ahead spot price (also known as the market clearing price), shortfalls are covered at an up-regulation price that is typically higher, and surpluses are cleared at a down-regulation price that is typically lower. Our scenario considers the standard European setting with two imbalance prices, under a price-taker assumption where individual wind producers do not influence market prices.

The nomination problem is complicated not only by the uncertainty of wind generation but also by the volatility of balancing prices. Up-regulation prices, in particular, can experience sudden spikes that lead to losses of exceptional magnitude. A large body of work has focused on forecasting production \cite{pinson2013wind, zhang2014review} and prices \cite{weron2014electricity, hong2016probabilistic}, and such models can be effective when extreme events exhibit early signals \cite{Christensen2012, galarneau2023foreseeing, sandhu2016forecasting}. However, many of the largest spikes arise from exogenous shocks such as generator outages or transmission faults, which are inherently unpredictable. These episodes meet all the criteria of \emph{black swans}: they are rare, unforeseeable, and highly consequential \cite{taleb2010black}. In our Finnish dataset, for instance, a single spring day in 2020 produced losses exceeding two weeks' profit within a single hour under a naive mean-forecast strategy.

Faced with this risk, a natural response is to adopt robust optimization, which guarantees performance even under the worst possible conditions \cite{ben2009robust}. While such strategies provide protection during extreme price spikes, they often achieve this by sacrificing profitability in normal hours. The resulting policies are safe but chronically conservative, earning less most of the time because the adverse events they guard against are exceedingly rare. A more effective hedge should therefore account not only for how severe deviations can be, but also for how likely they are. In other words, robustness must be calibrated to both magnitude and probability.

In this paper, we argue that \emph{optimal transport-based distributionally robust optimization} (OT-DRO) \cite{mohajerin2018data, blanchet2019quantifying, gao2023distributionally, shafiee2023nash} provides exactly the right balance between protection and performance. The key idea is to make decisions that remain effective even if the underlying probability distribution shifts within a small region defined by an optimal transport metric. This region captures both how far scenarios may deviate from historical data and how plausible those deviations are. Its two design parameters\,---\,the transportation cost and the radius\,---\,can be selected intuitively from data, and the resulting optimization problem admits a tractable convex reformulation. Using four years of Finnish wind and price data together with production forecasts from the Nordic energy company Fortum, we show that this approach consistently mitigates drawdowns during price spikes while remaining competitive in normal conditions. To the best of our knowledge, this is the first study to address price-spike hedging in day-ahead wind nomination using OT-DRO.

The rest of the paper is organized as follows. Section 2 formulates the nomination problem and quantifies the financial risks associated with price spikes using four years of market data. Section 3 develops the proposed OT-DRO framework, explains how its key parameters can be identified from historical data, and presents its convex reformulation. Section 4 reports the experimental validation, including case studies of extreme events and aggregate performance statistics, showing that robustification enhances resilience to black swans without compromising day-to-day profitability.



\section{Problem description}
\label{sec:description}

\begin{figure}[t]
\centering
\includegraphics[width=0.96\linewidth]{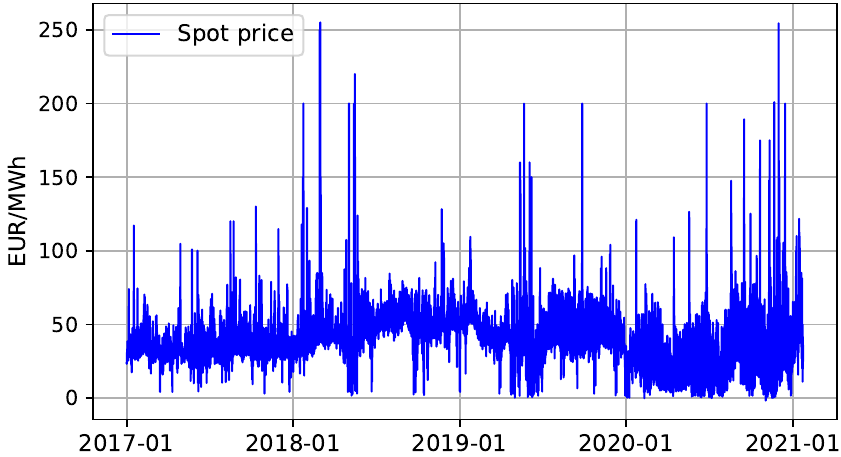}
\caption{Spot price from January 2017 to January 2021.}
\label{fig:spot:4years}
\end{figure}

\begin{figure}[t]
\centering
\includegraphics[width=\linewidth]{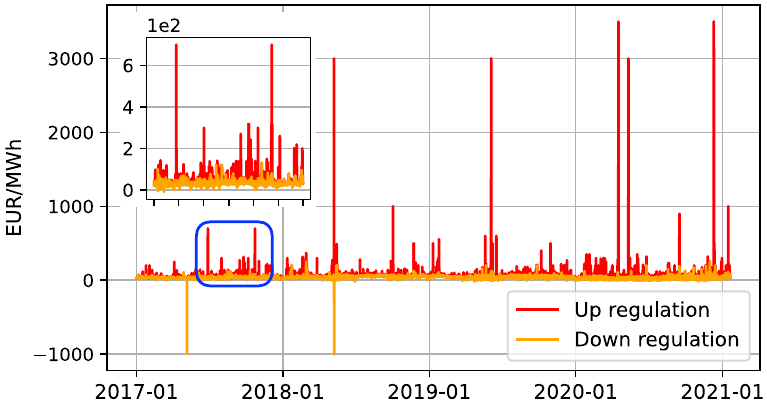}
\caption{Up- and down-regulation prices from January 2017 to January 2021.}
\label{fig:regulation:4years}
\end{figure}

Wind producers participating in the Nordic day-ahead market must nominate, by 12:00 (noon) each day, the amount of energy they will deliver during the 24-hour period of the following day. This results in 24 hourly nominations. To illustrate the mechanism, we focus on a single delivery hour and denote by $n$, $g$, $s$, $r^+$, and $r^-$ the nomination, realized generation, spot price, and the up- and down-regulation prices, respectively. The nomination is sold at the day-ahead spot price $s$, while real-time imbalances are settled at the regulation prices. If generation falls short of the nomination, the producer must buy the deficit at $r^+$; if generation exceeds the nomination, the excess is sold at $r^-$. These prices are determined ex post, after the delivery hour. The producer’s profit for that hour can therefore be written as
\begin{equation}
\pi(n, (g, s, r^-, r^+)) = n s + (g - n)_+ r^- - (n - g)_+ r^+ ,
\label{eq:profit}
\end{equation}
where $(x)_+ = \max(x,0)$. All price and generation data used in this study are obtained from Fingrid’s open data portal (\url{https://data.fingrid.fi/en/}), which provides hourly records for Finland. The evolution of the spot prices and of the up- and down-regulation prices from January 2017 to January 2021 is shown in Fig.~\ref{fig:spot:4years} and Fig.~\ref{fig:regulation:4years}, respectively. The general price relationship
\[
r^- < s < r^+
\]
holds for all hours. However, the regulation prices exhibit strong variability around this pattern. Fig.~\ref{fig:regulation:4years} also includes a zoomed-in view of a typical period without major spikes, illustrating that smaller fluctuations are present even in otherwise calm intervals. Occasionally, the regulation prices deviate drastically. The down-regulation price $r^-$ can even become negative, implying that a producer may have to pay for overproduction. In our four-year dataset, two negative peaks exceed 1000~€/MWh. Conversely, the up-regulation price $r^+$ exhibits several extreme spikes, with five exceeding 3000~€/MWh and six above 1000~€/MWh in total. Even smaller spikes are impactful: there are 236 hours when $r^+ > 200$~€/MWh, as visible in the zoomed portion of Fig.~\ref{fig:regulation:4years}, showing that moderate yet frequent surges also contribute significantly to financial risk.

To determine their day-ahead nominations, wind producers rely on generation forecasts that translate meteorological uncertainty into expected power output. In our case study, the company Fortum provided 52 ensemble forecasts for each hour, obtained by propagating different weather scenarios through a wind-generation model. Each trajectory represents a plausible realization of production, forming an empirical distribution of possible outputs.

In principle, additional forecasts could also be available for regulation prices. Such information could, in theory, support a more risk-averse nomination strategy by indicating when market volatility is expected to be high. However, the extreme price spikes that motivate this study are black swan events, triggered by exogenous contingencies such as grid faults or sudden outages, and are therefore inherently unpredictable. Even sophisticated forecasting systems cannot identify when such episodes will occur or how severe they will be.

One possible nomination strategy is to select the mean forecast, that is, the average of the ensemble predictions for each hour. This choice is natural when the forecasts are calibrated so that their mean matches the expected generation, and it serves as a useful benchmark for normal operating conditions. In practice, some producers adopt more risk-averse strategies that explicitly hedge against forecast errors or price spikes; however, it is unclear what a \emph{principled strategy} should be when confronted with black swan events. We return to this question and discuss robust approaches in Section~\ref{sec:methodology}. The mean-forecast policy nonetheless remains an informative baseline for understanding how extreme price spikes amplify financial risk. Since the deficit $(n - g)_+$ is multiplied by the up-regulation price $r^+$ in \eqref{eq:profit}, even modest forecast errors can translate into large losses during spike hours. Throughout the paper, we zoom in on spring 2020, a representative period containing two major price spikes. This quarter offers a detailed view of how extreme up-regulation episodes translate into large profit drawdowns and serves as the main empirical reference for our analysis.

\begin{figure}[t]
\centering
\includegraphics[width=\linewidth]{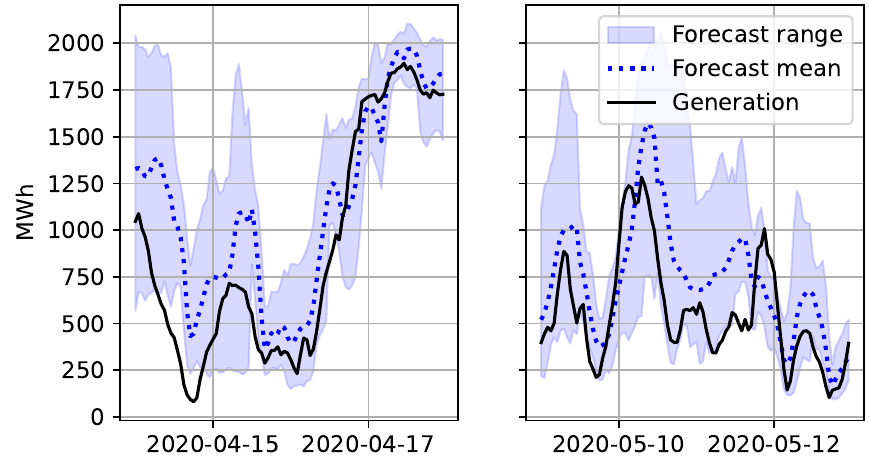}
\caption{Forecast range and actual generation during crash periods.}
\label{fig:forecast:crashes}
\end{figure}

\begin{figure}[t]
\centering
\includegraphics[width=\linewidth]{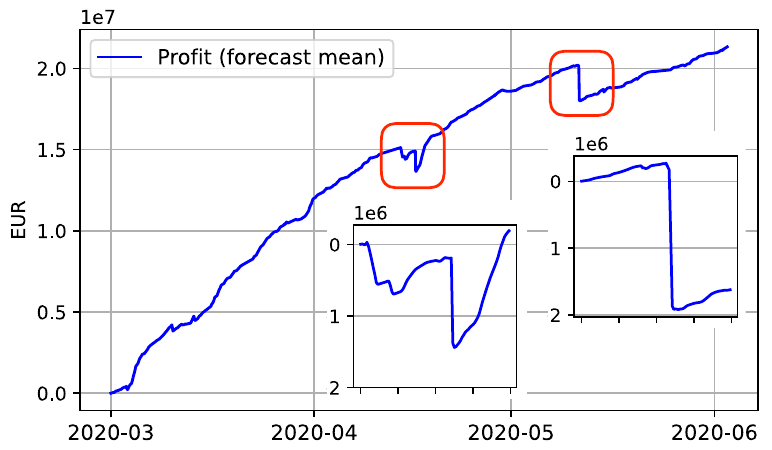}
\caption{Profit of the mean-forecast policy in spring 2020.}
\label{fig:profit:forecast}
\end{figure}

Fig.~\ref{fig:forecast:crashes} shows the range of the 52 forecasts around the two most extreme events of spring 2020, together with the realized generation. While the forecasts track the overall production pattern, they cannot anticipate the sudden up-regulation price spikes that turn routine forecast errors into large profit drops. The resulting profit, displayed in Fig.~\ref{fig:profit:forecast}, reveals two dramatic losses under the mean-forecast policy, the second exceeding €2 million within only a few hours. These crashes effectively wipe out several weeks of cumulative profits almost instantaneously, underscoring the disproportionate impact of black swans. Addressing this vulnerability is therefore essential not only for maintaining profitable operations but also for ensuring that renewable producers remain financially incentivized to participate in day-ahead markets despite the inherent volatility of wind generation.

\begin{figure}[t]
\centering
\includegraphics[width=\linewidth]{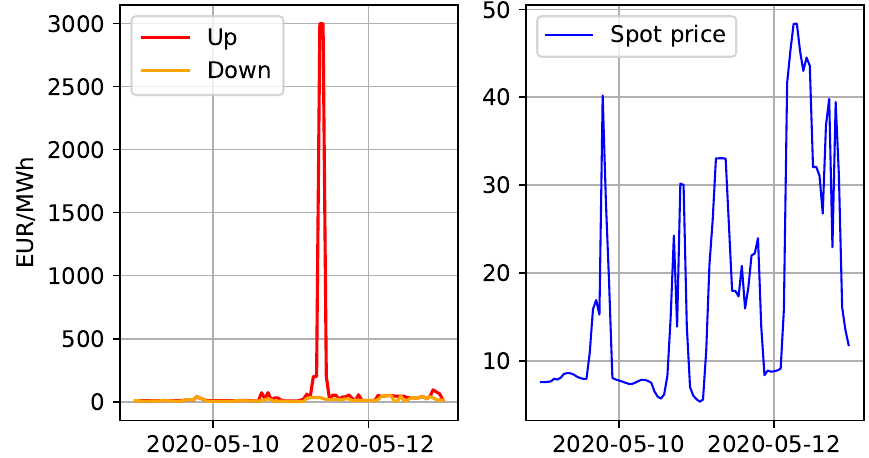}
\caption{Up- and down-regulation prices during the second crash.}
\label{fig:second:crash}
\end{figure}

Fig.~\ref{fig:second:crash} zooms in on the second major event of spring 2020, displaying the up- and down-regulation prices together with the spot price. The up-regulation price exceeds 3000~€/MWh during this episode, driving the largest loss observed under the mean-forecast policy. The plot highlights how abrupt and isolated these spikes are, reinforcing that such black swan events are virtually impossible to anticipate and dominate the overall profit distribution.


\section{Methodology}
\label{sec:methodology}

Since anticipating black swan events is impossible, the natural alternative is to be robust to them. Classical robust optimization aims to guarantee performance under the worst-case realization of uncertainty, but it implicitly treats all adverse outcomes as equally likely. This is problematic in energy markets, where the extreme price spikes observed in the data are rare. Designing an effective strategy therefore requires recognizing this imbalance: most fluctuations are mild, but a few can be catastrophic. A sound approach must balance magnitude and likelihood\,---\,a principle we refer to as \emph{robustness in both (space, likelihood)}.

Note that standard \emph{stochastic optimization} \cite{shapiro2021lectures} cannot address this challenge either, as it employs a single fixed probability distribution. While such a model can assign different likelihoods to adverse outcomes, these probabilities are fixed and cannot adapt to misspecification or unexpected events. In what follows, we argue that a \emph{distributionally robust} formulation\,---\,which can be viewed as an interpolation between robust and stochastic optimization\,---\,based on optimal transport can model precisely this flexibility, enabling robustness with respect to both spatial deviations and likelihood distortions.

The uncertainty in the day-ahead profit~\eqref{eq:profit} arises from four stochastic variables: the spot price, the up- and down-regulation prices, and the generation. We denote them collectively by $x = (g, s, r^-, r^+)$. The possible realizations of $x$ are modeled as a bounded set $\mathcal X \subset \mathbb{R}^4$, informed by historical data and expert judgment, representing feasible market conditions. Although the up-regulation price $r^+$ can reach high values, the Finnish data show it never exceeded 4000~€/MWh, suggesting a natural upper bound.

The joint uncertainty in $x$ can be characterized empirically from past data and forecasts. To make this precise, we construct an empirical distribution ${\mathbb P}_x$ over $\mathcal X$ that summarizes normal operating conditions. Each sample of ${\mathbb P}_x$ corresponds to one historical or forecasted realization of prices and generation, capturing their typical co-fluctuations. This empirical model reflects the variability observed in everyday market conditions but, by construction, cannot capture black swan behavior. The rare and extreme deviations seen in the data lie outside the support of ${\mathbb P}_x$; to reason about them, we will need a framework that quantifies how much and in what way the true distribution may differ from this empirical baseline.

\emph{Optimal transport} (OT) provides precisely such a framework: it quantifies how one probability distribution can differ from another by assigning a ``cost'' to moving probability mass between outcomes. It does so through two main ingredients. The first is a transportation cost $c(x_1, x_2)$, which specifies how expensive it is to shift probability mass from scenario $x_1$ to $x_2$; this cost encodes the geometry of plausible deviations. The second is a budget parameter $\varepsilon$, which limits the total amount of mass that can be displaced; this defines how much deviation from the empirical distribution is tolerated. Together, these two quantities determine the level of robustness to black swan events: the transportation cost controls how tails decay (\emph{robustness in likelihood}), while the budget $\varepsilon$ controls how far the distribution can move (\emph{robustness in space}).

\begin{figure}[t]
\centering
\includegraphics[width=0.75\linewidth]{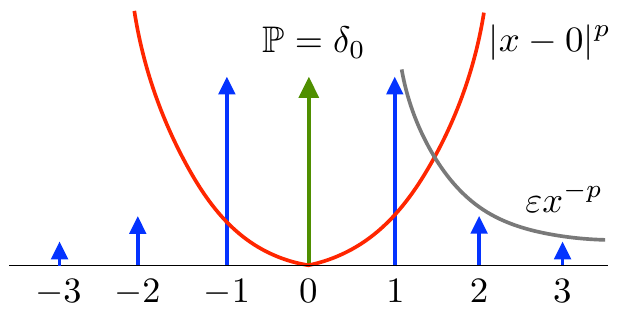}
\caption{Illustration of the maximum probability mass displacement from ${\mathbb P}_x$. The green arrow denotes the nominal distribution ${\mathbb P}_x$ supported only on the point $0$, the red curve represents the transportation cost $|\cdot|^p$, while the blue arrows denote the maximum probability mass that can be transported from ${\mathbb P}_x$ to that specific location. The grey curve denotes the upper bound on the probability displacements, inversely proportional to the transportation cost.}
\label{fig:displacement}
\end{figure}

Formally, the OT discrepancy between two probability distributions ${\mathbb P}_x, {\mathbb Q}_x \in \mathcal P(\mathcal X)$ with transportation cost function $c:\mathcal X\times\mathcal X\to\mathbb R_{\ge0}$ is defined as
\begin{align}
\label{eq:Wasserstein:distance}
    W({\mathbb P}_x,{\mathbb Q}_x)
    :=
    \inf_{\gamma\in\Gamma({\mathbb P}_x,{\mathbb Q}_x)} 
    \int_{\mathcal X \times \mathcal X} c(\xi_1, \xi_2)\, \text{d}\gamma(\xi_1,\xi_2),
\end{align}
where $\Gamma({\mathbb P}_x,{\mathbb Q}_x)$ is the set of all couplings with marginals ${\mathbb P}_x$ and ${\mathbb Q}_x$.
Despite this abstract definition, the idea is intuitive: $W({\mathbb P}_x,{\mathbb Q}_x)$ measures the minimum total cost of transporting the probability mass of ${\mathbb P}_x$ onto that of ${\mathbb Q}_x$.

In this framework, black swan behavior can be interpreted as transporting small probability mass to large deviations. Since the true probability of such extreme scenarios is unknown, it is natural to constrain it indirectly by limiting how much total transport can occur. This gives rise to the ambiguity set
\begin{align*}
    \mathbb B_\varepsilon({\mathbb P}_x)
    :=
    \left\{{\mathbb Q}_x\in\mathcal P(\mathcal X): W({\mathbb P}_x, {\mathbb Q}_x)\leq\varepsilon \right\},
\end{align*}
which collects all distributions that can be obtained from ${\mathbb P}_x$ by redistributing mass within a transport budget $\varepsilon$.
The radius $\varepsilon$ controls how much probability can move toward extreme events, while the choice of cost $c$ determines how rapidly that probability decays with distance, typically following a power-law behavior when $c(\xi_1,\xi_2)=\|\xi_1-\xi_2\|^p$, with $p \geq 1$.

To see this mechanism more concretely, consider the one-dimensional case with ${\mathbb P}_x = \delta_0$, that is $\xi_1 = 0$ with probability 1, and $c(\xi_1,\xi_2)=|\xi_1-\xi_2|^p$. Then, given any distribution ${\mathbb Q}_x$, there exists a unique coupling between $\delta_{0}$ and ${\mathbb Q}_x$, namely the product distribution $\delta_0 \otimes {\mathbb Q}_x$. Consequently, the minimum in \eqref{eq:Wasserstein:distance} drops and we recover $W(\delta_0, {\mathbb Q}_x) = \mathbb E_{\xi_2\sim{\mathbb Q}_x}[|\xi_2|^p]\leq\varepsilon$, limiting how much probability mass can be placed far from zero. If a fraction $\alpha$ of mass is shifted from $0$ to $x$, then $W({\mathbb P}_x,{\mathbb Q}_x)=\alpha|x|^p\leq\varepsilon$, implying $\alpha\leq\varepsilon |x|^{-p}$. Hence, the maximum probability assignable to a deviation of magnitude $x$ decays as $|x|^{-p}$ with prefactor $\varepsilon$. This decay mirrors a Pareto-type tail often used to describe black swan statistics: large shocks remain possible but become increasingly unlikely. Fig.~\ref{fig:displacement} illustrates this relationship. 

For instance, when $p=1$ and $\varepsilon=1$, at most $1/3000$ (or $0.03\%$) of the total probability mass can be placed at a deviation of $3000$, roughly matching the magnitude of the largest observed price spikes in the Finnish up-regulation data. Remarkably, the empirical decay of spike probabilities in Fig.~\ref{fig:regulation:4years} follows approximately an $x^{-1}$ pattern, consistent with this theoretical tail behavior. We will revisit this correspondence in detail in Section~\ref{sec:experiments}.

\begin{figure}[t]
    \centering
    \includegraphics[width=\linewidth]{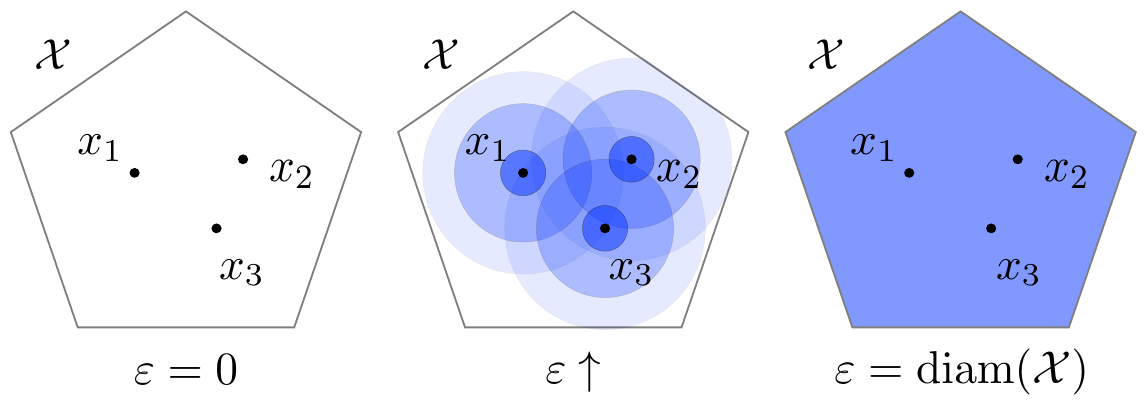}
    \caption{Illustration of the OT ambiguity set $\mathbb B_\varepsilon({\mathbb P}_x)$, with ${\mathbb P}_x = \frac{1}{3}\sum_{i=1}^3 \delta_{x_i}$, for increasing robustness radius. Each empirical sample $x_i$ can transport a limited amount of probability mass within the total budget, determined jointly by the transportation cost~$c$ and the radius~$\varepsilon$. Larger $\varepsilon$ values correspond to greater robustness, allowing more probability mass to move toward extreme scenarios. When $\varepsilon$ reaches the diameter of $\mathcal X$ (measured under the cost $c$), every scenario in $\mathcal X$ becomes attainable, recovering the classical robust optimization limit.}
    \label{fig:interpolation}
\end{figure}

This reasoning extends directly to empirical distributions supported on multiple samples, where each scenario transports a limited share of probability mass within the total budget $\varepsilon$, causing the ambiguity set to expand progressively from the empirical points toward the full support $\mathcal X$ (see Fig.~\ref{fig:interpolation}).

Before proceeding, we note that alternative choices exist for defining neighborhoods of probability distributions. One could, in principle, use other statistical distances such as the Kullback–Leibler (KL) divergence, total variation (TV) distance, or maximum mean discrepancy (MMD), to name just a few representative examples. However, these alternatives are less suited to capturing the phenomena described above. For example, KL divergence requires absolute continuity between distributions, which precludes modeling black swan events lying outside the support of ${\mathbb P}_x$ that represents normal operating conditions. Although TV and MMD do not suffer from this limitation, they lack the geometric and probabilistic interpretability of OT. Compared to the TV distance\,---\,which allows any fraction of probability mass to be reassigned arbitrarily\,---\,OT introduces a spatial penalty through the transportation cost, ensuring that distant deviations receive proportionally less mass. This spatial decay is also observed empirically in real markets; for instance, the probability of large price deviations in the Finnish up-regulation data (Fig.~\ref{fig:regulation:4years}) decreases roughly with distance, consistent with this modeling principle. MMD, in turn, measures discrepancies in a reproducing kernel Hilbert space. While theoretically appealing, its behavior depends intricately on the kernel choice, obscuring how probability mass shifts across scenarios. Moreover, as we will see in Section~\ref{sec:experiments}, OT leads to a computationally advantageous formulation of decision-making: the resulting optimization problem is convex\,---\,and in our case, even linear\,---\,whereas MMD-based robustness typically yields nonconvex problems.

Armed with this understanding of OT-based ambiguity sets and their ability to capture rare events, we now turn to the producer’s \emph{hourly} decision-making problem: determining a nomination strategy $n$ that remains profitable under distributional uncertainty. We formalize it as an \emph{OT-based distributionally robust optimization} (OT-DRO):
\begin{align}
\label{eq:dro}
    \max_{n} \; \min_{{\mathbb Q}_x \in \mathbb B_\varepsilon({\mathbb P}_x)} \mathbb E_{x \sim {\mathbb Q}_x} \left[\pi(n, x)\right],
\end{align}
where $\pi(n, x)$ denotes the profit defined in~\eqref{eq:profit}, and $\mathbb B_\varepsilon({\mathbb P}_x)$ represents the OT ambiguity set around the empirical distribution ${\mathbb P}_x$. Problem~\eqref{eq:dro} can be interpreted as a two-player, zero-sum game between the producer and nature. The producer selects the nomination $n$ to maximize profit, while nature, acting adversarially, perturbs the empirical distribution within the transportation budget $\varepsilon$ to minimize it. Small values of $\varepsilon$ correspond to optimistic, data-driven strategies that rely heavily on the empirical distribution, whereas larger values reflect increasing robustness to black swan shifts.

Despite the infinite-dimensional minimax structure of the OT-DRO problem~\eqref{eq:dro}, it is known that such problems admit exact convex reformulations under mild regularity conditions on the cost function and the loss. In particular, \cite{mohajerin2018data} shows that when the transportation cost is induced by a norm (i.e., $c(x_1,x_2)=\|x_1-x_2\|$), the inner problem admits a strong dual representation, leading to a tractable convex program. More recently, \cite{shafiee2023nash} extended this computational theory to arbitrary convex costs (encompassing $c(x_1,x_2)=\|x_1-x_2\|^p$, for any $p \geq 1$, thereby covering the class of power-law behaviors that connect naturally to the Pareto-type decay discussed above). These reformulations make OT-based distributional robustness both conceptually and computationally viable for real-time decision-making in energy markets. In the next section, we calibrate the OT model to Finnish market data, showing how empirical spike behavior guides the selection of both the cost function and robustness level.


\section{Experiments}
\label{sec:experiments}

We now provide a detailed empirical analysis of how to select the key ingredients of the OT-DRO model: 
\begin{itemize}
    \item the reference distribution~${\mathbb P}_x$, 
    \item the support set $\mathcal X$,
    \item the transportation cost~$c$, and 
    \item the robustness radius~$\varepsilon$.
\end{itemize}
To evaluate the robustness of our method under realistic conditions, we divide the four-year Finnish dataset into seasons, each spanning three consecutive months, and treat each season once as the test period while using the remaining data as training. Using this setup, we show how each component of the model can be grounded in observed market behavior, leading to a data-driven formulation of the robust nomination problem.

\subsection{Reference distribution} We aim to construct an empirical distribution ${\mathbb P}_x$ over $x=(g,s,r^-,r^+)$ that reflects realistic market conditions given the producer’s available information. While producers could, in principle, forecast the full joint distribution of these variables, in our case only a forecast $f$ for generation $g$ is available. We use this to infer a plausible empirical distribution for $x$ by conditioning on similar past forecasts. Formally, let $\{(f_i,g_i,s_i,r_i^-,r_i^+)\}_{i=1}^N$ denote the training data. Given a new forecast $f$, we select the $\alpha$ fraction of historical samples whose forecasts $f_i$ are closest to $f$ in absolute difference, reasoning that similar forecasts correspond to comparable operating conditions. Denote by $\mathcal I$ the set of selected indices, with $|\mathcal I| = \lceil \alpha N \rceil$. The empirical distribution is then defined as
\begin{equation*}
    {\mathbb P}_x = \sum_{i\in\mathcal I} w_i \, \delta_{(g_i,s_i,r_i^-,r_i^+)},
\end{equation*}
where the weights are chosen according to
\begin{equation*}
    w_i \propto \left(1 - \frac{\|f-f_i\|}{d_{\max}}\right)^{\beta},
\end{equation*}
and $d_{\max}$ is the largest distance among the selected neighbors. The parameter $\beta>0$ controls how sharply the weights decay with forecast distance. In our experiments, we use $\alpha=1/3$ and $\beta=2$. This construction parallels a nearest-neighbor conditional density estimator, producing a data-driven distribution consistent with both historical variability and current forecasts.

\begin{table}[t]
\centering
\caption{Empirical exceedances of the up-regulation price $r^+$ over four years of Finnish market data. The probabilities of large deviations decay approximately as $q^{-1}$, consistent with a Pareto-type tail corresponding to $p=1$ in the OT model. The implied $\varepsilon$ values justify the three robustness levels $\varepsilon\!\in\!\{0.5,1.0,1.5\}$ used in our experiments.}
\setlength{\tabcolsep}{8pt}
\renewcommand{\arraystretch}{1.2}
\begin{tabular}{|c|c|c|c|}
\hline
\textbf{$q$ (€/MWh)} & \textbf{Count $r^+\!\ge q$} & \textbf{Freq. (\%)} & \textbf{$\text{Pr}\{r^+\!\ge q\} \approx \varepsilon q^{-1}$} \\
\hline
200  & 236 & 0.68 & $1.3\,q^{-1}$ \\
\hline
500  & 20  & 0.06 & $0.3\,q^{-1}$ \\
\hline
1000 & 7   & 0.02 & $0.2\,q^{-1}$ \\
\hline
3000 & 5   & 0.01 & $0.4\,q^{-1}$ \\
\hline
\end{tabular}
\label{tab:tail_rplus}
\end{table}

\subsection{Support set} 

We define the uncertainty support $\mathcal X$ as a bounded box encompassing all plausible realizations of $x=(g,s,r^-,r^+)$. Each component is delimited by its historical range, enlarged by a 20\% safety margin to accommodate potential variability. Specifically, the variables $g$, $s$, and $r^+$ are allowed to vary between 0.8 times their historical minimum and 1.2 times their historical maximum, while for the down-regulation price $r^-$, which can take negative values, the lower bound is set to 1.2 times its minimum. This construction ensures that $\mathcal X$ captures realistic yet sufficiently broad market conditions without extrapolating to implausible extremes. While bounding $\mathcal X$ is not theoretically necessary\,---\,none of our results depend on this restriction\,---\,it prevents overly conservative solutions that could arise from assigning mass to completely implausible events.

\begin{figure}[t]
    \centering
    \includegraphics[width=\linewidth]{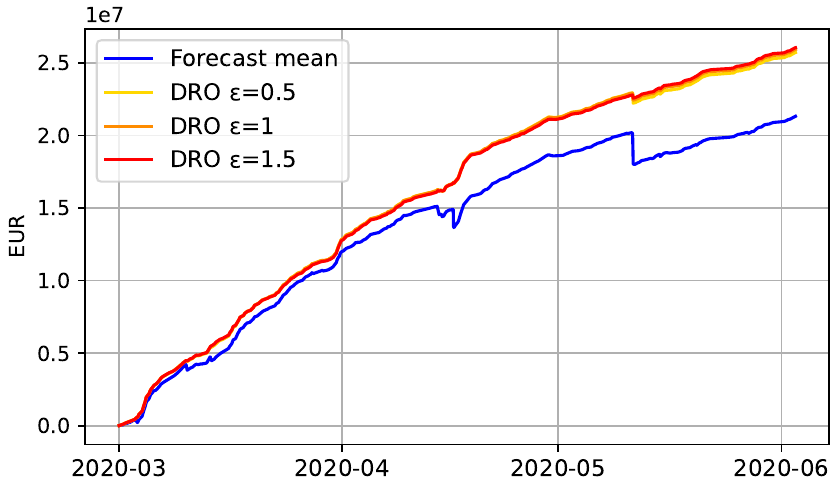}
    \caption{Comparison of the profit trajectories for the forecast-mean strategy and the OT-DRO strategies with radii $\varepsilon\in\{0.5,1,1.5\}$ during spring~2020. Larger robustness radii mitigate the impact of the two market crashes, confirming the protective effect of higher $\varepsilon$.}
    \label{fig:profit:DRO_profit}
\end{figure}

\subsection{Transportation cost and radius}

We adopt a transportation cost of the form $c(\xi_1,\xi_2)=\|\xi_1-\xi_2\|^p$, for some norm (e.g., $\ell_1$ or $\ell_\infty$), where the exponent $p$ controls how strongly distant deviations are penalized. To identify an appropriate value for~$p$, we focus on the up-regulation price $r^+$, whose extreme variability dominates market risk. We can therefore restrict the reasoning to one dimension. Let $q$ denote a realization of $r^+$. As shown in Section~\ref{sec:methodology}, the OT constraint implies that a fraction $\alpha$ of probability mass placed at deviation $q$ must satisfy $\alpha q^p \le \varepsilon$, or equivalently $\alpha \le \varepsilon q^{-p}$. Thus, the parameter $p$ determines how fast the assignable probability decays with the magnitude of deviation. 

Examining the empirical tail of $r^+$ over roughly $35{,}000$ hourly observations, we find that 236 prices exceed 200~€/MWh, 20 exceed 500~€/MWh, 7 exceed 1000~€/MWh, and 5 exceed 3000~€/MWh. The corresponding exceedance probabilities\,---\,$0.68\%$, $0.06\%$, $0.02\%$, and $0.01\%$\,---\,decrease approximately as $q^{-1}$, consistent with a Pareto-like envelope for $p=1$. This motivates adopting $c(\xi_1,\xi_2)=\|\xi_1-\xi_2\|$ as the transportation cost, reproducing the empirically observed heavy-tailed behavior of $r^+$. Moreover, using the relation $\alpha \le \varepsilon q^{-1}$, we can interpret $\varepsilon$ as the effective transport budget required to sustain an empirical tail probability $\alpha$ at deviation $q$. Identifying $\alpha$ with the observed exceedance $\text{Pr}\{r^+\!\ge q\}$ yields $\varepsilon \approx q\,\text{Pr}\{r^+\!\ge q\}$. Applying this relation to the Finnish data gives approximate values $\varepsilon_{200}\!\approx\!1.3$, $\varepsilon_{500}\!\approx\!0.3$, $\varepsilon_{1000}\!\approx\!0.2$, and $\varepsilon_{3000}\!\approx\!0.4$ (see Table~\ref{tab:tail_rplus}). Motivated by this, in our experiments we focus on three robustness radii, $\varepsilon\in\{0.5,1,1.5\}$, corresponding to realistic low, medium, and high levels of robustness.

\begin{figure}[t]
    \centering
    \includegraphics[width=\linewidth]{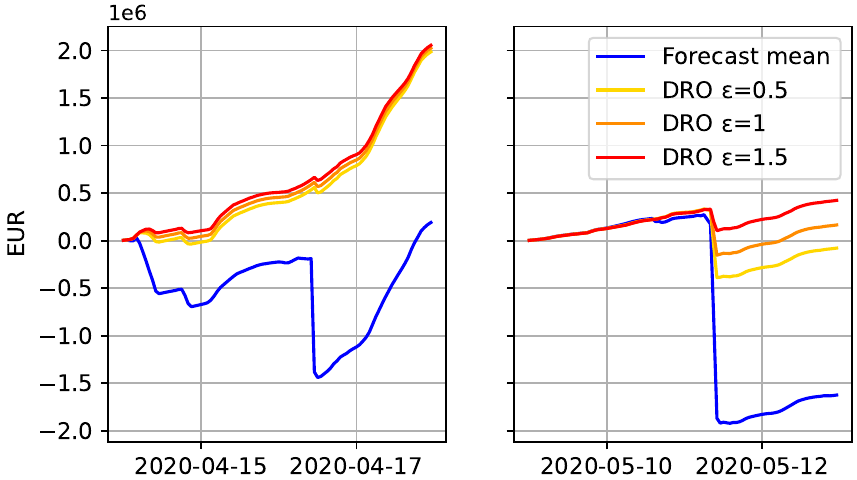}
    \caption{Zoom-in on the two major crashes of spring~2020, comparing the profit trajectories of the forecast-mean and OT-DRO strategies. Larger radii $\varepsilon$ provide smoother and less severe drawdowns, highlighting the improved robustness of the DRO approach during extreme market conditions.}
    \label{fig:profit:DRO_crashes}
\end{figure}

\subsection{Convex reformulation}
\label{sec:finite-convex}

The OT-DRO problem in~\eqref{eq:dro} is infinite-dimensional, as the inner minimization ranges over all probability distributions within the OT ambiguity set. Nonetheless, following the results of~\cite{mohajerin2018data, shafiee2023nash}, such problems can be reformulated exactly as finite-dimensional convex programs when the objective function is convex or concave piecewise affine in the uncertain variables. Observing that the profit can be equivalently rewritten as
\[
\pi\big(n,(g,s,r^-,r^+)\big) = n s - \max\{(n-g)r^+, (n-g)r^-\},
\]
we expand the uncertainty vector to $\xi := (g,s,r^-,r^+,t^-,t^+)$ by introducing the auxiliary variables $t^- := g r^-$ and $t^+ := g r^+$. This representation makes the profit a concave piecewise linear function of~$\xi$. Because this change of variables is a deterministic transformation of the original uncertainty, the key ingredients of the OT-DRO formulation can be recovered directly from their counterparts in the $(g,s,r^-, r^+)$ coordinates, as detailed below:

\smallskip

\noindent\emph{Reference distribution.} 
The empirical reference $\widehat{\mathbb P}_\xi$ is obtained directly from $\widehat{\mathbb P}_x$ by lifting each sample through the deterministic map $\Phi(x)=(x,\,g r^-,\,g r^+)$. If $(g_i,s_i,r_i^-,r_i^+)$ appears in $\widehat{\mathbb P}_x$ with weight $w_i$, then $(g_i,s_i,r_i^-,r_i^+,\,g_i r_i^-,\,g_i r_i^+)$ appears in $\widehat{\mathbb P}_\xi$ with the same weight. The two empirical distributions therefore share the same samples and probabilities, differing only by the inclusion of the auxiliary coordinates $t^-$ and $t^+$.

\smallskip

\noindent\emph{Support set.}
The support $\Xi$ is defined analogously to $\mathcal X$ as a box containing all plausible realizations of $\xi$. Each coordinate is bounded by its historical range enlarged by a 20\% margin, with the same enlargement applied to the empirical ranges of the products $g r^-$ and $g r^+$. Although bounding $\Xi$ is not required by theory, it provides a practical safeguard against allocating probability mass to implausible extremes. Formally, we represent this box as a polyhedron
\[
\Xi:= \{\xi \in \mathbb R^6:\; C \xi \leq d\},
\]
where $C\in\mathbb R^{12\times6}$ and $d\in\mathbb R^{12}$ collect the upper and lower bounds on each coordinate of~$\xi$.

\smallskip

\noindent\emph{Transportation cost and radius.}
We retain the same transportation cost $c(\xi_1,\xi_2)=\|\xi_1-\xi_2\|$ and robustness radii $\varepsilon\in \{0.5,1,1.5\}$. The transformation $\Phi(x)=(x,gr^-,gr^+)$ is deterministic and Lipschitz on the bounded set $\mathcal X$, so it represents a reparametrization rather than an expansion of the uncertainty. Moreover, by introducing the auxiliary variables $t^\pm$ without enforcing the exact relations $t^\pm=g r^\pm$ within the support, the feasible set $\Xi$ becomes a slightly larger relaxation of the original uncertainty domain. For a fixed $c$ and $\varepsilon$, this enlargement can only make the formulation weakly more conservative, thus preserving or slightly increasing robustness. Keeping the same $\varepsilon$ therefore provides a consistent and safe level of protection against distributional perturbations without requiring any rescaling of the transport cost.

\smallskip

With these elements in place, the OT-DRO problem of interest can be equivalently written as
\begin{align}
\label{eq:dro:updated}
    \max_{n} \; \min_{{\mathbb Q}_\xi \in \mathbb B_\varepsilon({\mathbb P}_\xi)} \mathbb E_{\xi \sim {\mathbb Q}_\xi} \left[\pi(n, \xi)\right],
\end{align}
where the profit, expressed in the expanded coordinates $\xi=(g,s,r^-,r^+,t^-,t^+)$, takes the piecewise linear form
\begin{align*}
    \pi(n,\xi) &= -\max\{n r^+ - t^+ - n s, n r^- - t^- - n s\}\\
    &= -\max\{a_1(n)^\top \xi, a_2(n)^\top \xi\},
\end{align*}
with nomination-dependent $a_1(n):=(0,-n,0,n,0,-1)$ and $a_2(n):=(0,-n,n,0,-1,0)$. Since the transportation cost is induced by a norm, \cite[Corollary 5.1(i)]{mohajerin2018data} applies directly: taking the negative sign outside turns~\eqref{eq:dro:updated} into the standard minimax form, which admits an exact finite-dimensional convex reformulation yielding the optimal nomination~$n^\star$:
\begin{align*}
    \begin{array}{cll}
        n^\star = \underset{\substack{n, \lambda, s_i,\\ \gamma_{i1}, \gamma_{i2}}}{\text{arg}\min} & \lambda \varepsilon +\sum_{i=1}^{\lceil \alpha N \rceil} w_i s_i \\
        [3.5ex]~\text{ \quad }\st
        & \DS 
        a_1(n)^\top \xi_i + \gamma_{i 1}^\top (d - C \xi_i) \leq s_i &\forall i \leq \lceil \alpha N \rceil
        \\
        [1ex]
        & \DS
        a_2(n)^\top \xi_i + \gamma_{i 2}^\top (d - C \xi_i) \leq s_i &\forall i \leq \lceil \alpha N \rceil
        \\
        [1ex]
        & \DS 
        \left\| C^\top \gamma_{i 1} - a_1(n) \right\|_* \leq \lambda
        &\forall i \leq \lceil \alpha N \rceil
        \\
        [1ex]
        & \DS 
        \left\| C^\top \gamma_{i 2} - a_2(n) \right\|_* \leq \lambda
        &\forall i \leq \lceil \alpha N \rceil
        \\
        [1ex]
        &  \DS
         \gamma_{i1} \geq 0, \gamma_{i2} \geq 0
         &\forall i \leq \lceil \alpha N \rceil
    \end{array}
\end{align*}
with $\|\cdot\|_*$ being the dual norm of $\|\cdot\|$. In our experiments, we selected $\alpha = 1/3$ and $\|\cdot\|_1$ so that $\|\cdot\|_* = \|\cdot\|_\infty$. Note that with this choice of norm, the reformulation reduces to a \emph{linear program}.

\subsection{Empirical Performance}

Applying the reformulation to the Finnish dataset yields exactly the behavior anticipated by our modeling choices: the OT-DRO nominations hedge black swans while remaining competitive in normal hours. Fig.~\ref{fig:profit:DRO_profit} shows that the cumulative profit in spring~2020 under OT-DRO stays above the forecast-mean benchmark and, crucially, avoids the large drops during the two spike episodes. The zoomed view in Fig.~\ref{fig:profit:DRO_crashes} highlights this effect: as the robustness radius increases from $\varepsilon=0.5$ to $1$ and $1.5$, the trajectory becomes progressively more resilient, with markedly smaller drops at both events. This monotone trend corroborates our calibration: once the transportation cost matches the empirical tail of up-regulation prices, the radius $\varepsilon$ acts as the correct “knob’’ to tune robustness against black swans.

\begin{table}[t]
\centering
\caption{Seasonal profit difference (\%) of OT-DRO strategies with $\varepsilon \in \{0.5, 1, 1.5\}$ relative to the forecast-mean strategy. Positive values indicate higher profit for OT-DRO.}
\setlength{\tabcolsep}{8pt} 
\renewcommand{\arraystretch}{1.3} 
\begin{tabular}{|l|c|c|c|}
\hline
\textbf{Season} & $\varepsilon=0.5$ & $\varepsilon=1$ & $\varepsilon=1.5$ \\
\hline
Winter 2016-17 & $-0.53\%$ & $-0.86\%$ & $-1.24\%$ \\
\hline
Spring 2017     & $-1.68\%$ & $-2.06\%$ & $-2.48\%$ \\
\hline
Summer 2017     & $+0.82\%$  & $+0.84\%$  & $+0.72\%$  \\
\hline
Autumn 2017     & $-0.19\%$ & $-0.20\%$ & $-0.29\%$ \\
\hline
Winter 2017-18 & $-0.37\%$ & $-0.57\%$ & $-0.81\%$ \\
\hline
Spring 2018     & $-1.70\%$ & $-2.09\%$ & $-2.42\%$ \\
\hline
Summer 2018     & $-0.10\%$ & $-0.21\%$ & $-0.35\%$ \\
\hline
Autumn 2018     & $-0.12\%$ & $-0.20\%$ & $-0.30\%$ \\
\hline
Winter 2018-19 & $+0.15\%$  & $+0.05\%$  & $-0.10\%$ \\
\hline
Spring 2019     & $-0.16\%$ & $-0.33\%$ & $-0.54\%$ \\
\hline
Summer 2019     & $-0.32\%$ & $-0.54\%$ & $-0.83\%$ \\
\hline
Autumn 2019     & $+0.78\%$  & $+0.65\%$  & $+0.44\%$  \\
\hline
Winter 2019-20 & $-0.20\%$ & $-0.55\%$ & $-0.99\%$ \\
\hline
\textbf{Spring 2020} & {\color{red}\textbf{+21.04\%}} & {\color{red}\textbf{+22.06\%}} & {\color{red}\textbf{+22.64\%}} \\
\hline
Summer 2020     & $+0.60\%$  & $+0.30\%$  & $-0.18\%$ \\
\hline
Autumn 2020     & $+0.10\%$  & $-0.15\%$ & $-0.54\%$ \\
\hline
Winter 2020-21 & $-0.03\%$ & $+0.03\%$  & $-0.03\%$ \\
\hline
\end{tabular}
\label{tab:profit_dro}
\end{table}

However, a natural question arises from Fig.~\ref{fig:profit:DRO_profit}: \emph{why does the OT-DRO profit curve remain consistently above the forecast-mean benchmark, even outside the two major crashes?} At first glance, this seems counterintuitive\,---\,robust formulations are typically more conservative, not systematically more profitable. The explanation lies in the structure of spring~2020: beyond the two black swan events, the up-regulation prices exhibit a dense sequence of smaller spikes, with amplitudes around 200–300~€/MWh (see Fig.~\ref{fig:regulation:4years}). These moderate but frequent deviations make robustness advantageous throughout the entire season. By hedging against local price excursions as well as extreme ones, the OT-DRO strategy achieves higher stability and, in aggregate, higher realized profit\,---\,appearing as a ``free lunch'' that in fact reflects the persistent volatility of that period. 

Beyond the turbulent spring~2020, the OT-DRO nomination behaves as expected under normal operating conditions. Table~\ref{tab:profit_dro} reports the seasonal profit differences between OT-DRO and the forecast-mean strategy across all four years and for the three robustness radii. In well-behaved seasons\,---\,when price spikes are infrequent or moderate\,---\,the OT-DRO strategy is slightly more conservative, yielding profits marginally below those of the forecast-mean benchmark. This modest cost of robustness, however, is largely compensated during volatile periods: in spring~2020, OT-DRO achieves a $21$–$22\%$ profit increase, effectively neutralizing the impact of the two black swan events. While black swan episodes occurred in other seasons as well (a total of five major spikes over the four years), OT-DRO did not yield a significant advantage in all of them. The reason is that the model’s robustness is calibrated \emph{in distribution}, meaning it anticipates tail risks on average rather than guaranteeing protection against every individual extreme. Achieving such pathwise robustness would require a deterministic worst-case formulation\,---\,equivalent to setting $\varepsilon$ to the diameter of~$\Xi$ (on the order of thousands)\,---\,which would make the strategy prohibitively conservative and unprofitable. Instead, OT-DRO strikes a balanced compromise: it slightly sacrifices profit in calm markets while maintaining strong protection against the typical scale and frequency of black swan events observed in real data. 

\begin{figure}[t]
    \centering
    \includegraphics[width=\linewidth]{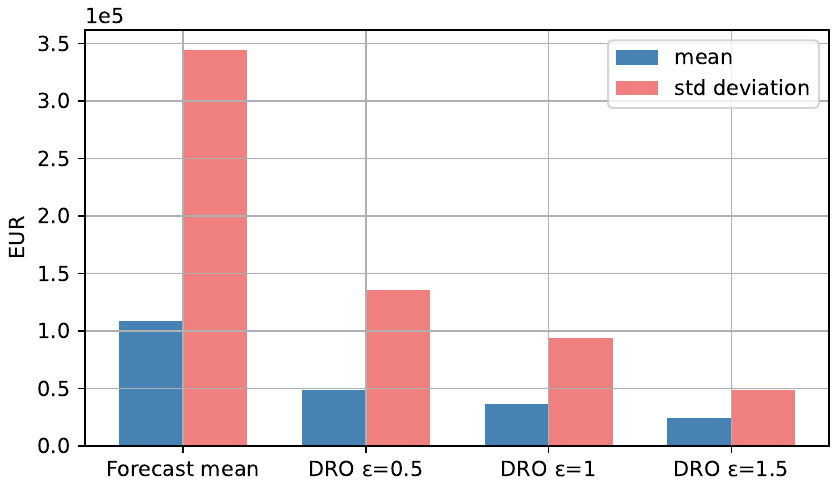}
    \caption{Summary statistics of profit drops over the four-year period. Each drop is defined from the onset of a decline in cumulative profit until the first subsequent recovery. The OT-DRO strategies yield smaller mean and standard deviation of losses compared to the forecast-mean strategy, confirming the stabilizing effect of robustness.}
    \label{fig:profit:DRO_drops}
\end{figure}

This behavior is further confirmed in Fig.~\ref{fig:profit:DRO_drops}, which summarizes the statistics of the profit drops for the three radii. Each drop is defined from the onset of a decline in cumulative profit until the first subsequent recovery. As the robustness radius~$\varepsilon$ increases, both the mean and standard deviation of the drops decrease, demonstrating that OT-DRO not only mitigates the largest losses but also stabilizes smaller fluctuations in profit.


\section{Conclusion}

This work introduced an OT-DRO formulation for day-ahead wind nominations that anticipates black swan events while maintaining profitability under normal conditions. Our formulation contributes to the growing line of research employing distributionally robust methods in power system operations and renewable trading \cite{Pinson2023, wang2018risk, poolla2020wasserstein, guo2018data}. However, unlike prior approaches that motivate robustness through forecast or model uncertainty, none have explicitly addressed protection against black swan events.

A promising extension is to adapt the approach to the single-price imbalance system recently adopted in several European markets, modeling a portfolio of spatially correlated wind farms that can exchange energy internally, and to extend the framework to the price-making regime\,---\,where learning becomes essential as the nomination itself influences market prices \cite{singhal2025learn}.

\bibliographystyle{unsrt}
\bibliography{references}
\end{document}